\input phyzzx
\year=1996
\catcode `\@=11
\newskip\frontpageskip
\newtoks\Ftuvnum   \let\ftuvnum=\Ftuvnum	
\newtoks\Ificnum   \let\ificnum=\Ificnum
\newtoks\Pubtype  \let\pubtype=\Pubtype
\newif\ifp@bblock  \p@bblocktrue
\def\PH@SR@V{\doubl@true \baselineskip=18.1pt plus 0.1pt minus 0.1pt
             \parskip= 3pt plus 2pt minus 1pt }
\def\PHYSREV{\papers\PhysRevtrue\PH@SR@V}
\let\physrev=\PHYSREV
\def\titlepage{\FRONTPAGE\papers\ifPhysRev\PH@SR@V\fi
   \ifp@bblock\p@bblock \else\hrule height\z@ \rel@x \fi }
\def\nopubblock{\p@bblockfalse}
\def\endpage{\vfil\break}
\frontpageskip=12pt plus .5fil minus 2pt
\Pubtype={}
\Ftuvnum={}
\Ificnum={}
\def\p@bblock{\begingroup \tabskip=\hsize minus \hsize
   \baselineskip=1.5\ht\strutbox \topspace-2\baselineskip
   \halign to\hsize{\strut ##\hfil\tabskip=0pt\crcr
       \the\Ftuvnum\crcr\the\Ificnum\crcr
	\the\date\crcr\the\pubtype\crcr}\endgroup}

\newcount\anni
\anni=\year
\advance\anni by -1900
\def\aapub{\afterassignment\aap@b\toks@}
\def\aap@b{\edef\n@xt{\Ftuvnum={FTUV\ \the\anni--\the\toks@}}\n@xt}

\def\bbpub{\afterassignment\bbp@b\toks@}
\def\bbp@b{\edef\n@xt{\Ificnum={IFIC\ \the\anni--\the\toks@}}\n@xt}

\let\ftuvnum=\aapub
\let\ificnum=\bbpub
\def\ftuvbin{
      \expandafter\ifx\csname binno\endcsname\relax%
         \expandafter\ifx\csname MailStop\endcsname\relax%
         \else%
            , Mail Stop \MailStop%
         \fi%
      \else%
         , Mail Stop \binno%
      \fi%
   }
\expandafter\ifx\csname eightrm\endcsname\relax
      \fi
\def\FTUV{\address{Departamento de F\'{\i}sica Te\'{o}rica and IFIC,\break
      Centro Mixto Univ. de Valencia-CSIC\break 
	46100-Burjassot (Valencia), Spain}}
\VOFFSET=33pt
\papersize
\physrev

\font\black=msbm10 scaled\magstep1

\def\NAME{\author{J.A. de Azc\'{a}rraga
\footnote{\tenrm\dag}{E-mail: azcarrag@evalvx.ific.uv.es}, 
A. M. Perelomov
\footnote{\tenrm\star}
{On leave of absence from Institute for Theoretical and
Experimental Physics, 117259 Moscow, Russia.\quad
E-mail: perelomo@evalvx.ific.uv.es} and 
J.C. P\'{e}rez Bueno \footnote{\tenrm\ddag}{E-mail: pbueno@lie.ific.uv.es}}}

\catcode `\@=\active

\def\ld{\mathop{\ldots}\limits}
\def\frac#1#2{{#1\over#2}}

\def\campo{{\cal X}}

\def\field #1{\hbox{{\black #1}}}

\def\pois#1#2{\{#1,#2\}}         
\def\dd#1{\frac{\partial}{\partial#1}}
\def\set#1{\{\,#1\,\}}         

\def\R{{\hbox{{\field R}}}}

\def\g{{\cal G}}

\REF\Na{Nambu, Y.: 
{\it Generalized Hamiltonian dynamics},
Phys. Rev. {\bf D7}, 2405--2412 (1973)}

\REF\BF{Bayen, F. and Flato, M.: 
{\it Remarks concerning Nambu's generalized mechanics}, 
Phys. Rev. {\bf D11}, 3049--3053 (1975)}

\REF\MS{Mukunda, N. and Sudarshan E.: 
{\it Relation between Nambu and Hamiltonian mechanics}, 
Phys. Rev. {\bf D13}, 2846--2850 (1976)}

\REF\Hir{Hirayama, H.: 
{\it Realization of Nambu mechanics: A particle 
interacting with an $SU(2)$ monopole}, 
Phys. Rev. {\bf D16}, 530--532 (1977)}

\REF\Ta{Takhtajan, L.:
{\it On foundations of the generalized Nambu mechanics},
Commun. Math. Phys. {\bf 160}, 295--315 (1994)}

\REF\Cha{Chatterjee, R.: {\it Dynamical symmetries and Nambu mechanics}, Stony 
Brook preprint (Jan. 1995)}

\REF\Lich{Lichnerowicz, A.:
{\it Les vari\'et\'es de Poisson et leurs alg\`ebres de Lie associ\'ees},
J. Diff. Geom. {\bf 12}, 253-300 (1977)} 

\REF\We{Weinstein, A.:
{\it The local structure of Poisson manifolds},
J. Diff. Geom. {\bf 18}, 523--557 (1983)}

\REF\BFFLS{Bayen, F.; Flato, M.; Fronsdal, C.; Lichnerowicz, A. 
and Sternheimer, D.: 
{\it Deformation theory and quantization},
Ann. Phys. {\bf 111}, 61--151 (1978)}

\REF\Li{Lie, S.:
{\it Begr\"undung einer Invariantentheorie der Ber\"uhrungs Transformationen},
Math. Ann. {\bf 8}, 214--303 (1874/75)}

\REF\Lie{Lie, S. and Engel, F.: 
{\it Theorie der Transformationsgruppen I--III}, 
Teubner (1888) (Chelsea, 1970)}

\REF\GMP{Grabowski, J.; Marmo, G. and Perelomov, A. M.:
{\it Poisson structures: towards a classification},
Mod. Phys. Lett. {\bf A18}, 1719--1733 (1993)}

\REF\CIMP{Cari\~ nena, J.; Ibort, A.; Marmo, G. and Perelomov A. M.:
{\it On the geometry of Lie algebras and Poisson tensors},
J. Phys. {\bf A27}, 7425--7449 (1994)}

\REF\AP{Alekseevsky, D.V. and Perelomov, A.M.:
{\it Poisson brackets on Lie algebras} 	   
Preprint ESI 247 (1995)}

\REF\Tu{Tulczyjev, W.M.: 
{\it Poisson brackets and canonical manifolds}, 
Bull. Acad. Pol. Sci. (Math. and Astronomy) {\bf 22}, 931--934 (1974)}

\REF\Sc{Schouten, J.A.:
{\it Ueber Differentialkomitanten zweir kontravarianter Gr\"oszen}, 
Proc. Kon. Ned. Akad. Wet. Amsterdam {\bf 43}, 449-452 (1940)}

\REF\Ni{Nijenhuis, A.: 
{\it Jacobi-type identities for bilinear differential concomitants 
of certain tensor fields}, 
Indag. Math. {\bf 17}, 390-403 (1955)}

\REF\APPB{de Azc\'arraga, J. A.; Perelomov, A. M. and P\'erez Bueno, J. C.:
in preparation}

\REF\Po{Poisson, S.:
{\it M\'emoire sur la variation des constantes arbitraires dans les 
questions de m\'echanique},
J. \'Ecole Polytec. {\bf 8}, 266--344 (1809)}

\REF\CE{Chevalley, C. and Eilenberg, S.: {\it Cohomology theory of
Lie groups and Lie algebras}, Trans. Am. Math. Soc. {\bf 63}, 85--124 
(1948)} 

\REF\AI{de Azc\'arraga, J. A. and Izquierdo, J.M.: {\em Lie algebras, Lie 
groups, cohomology and some applications in physics}, Camb. Univ. Press 
(1995)}

%
%
%
%
%
%
%
%
%

\ftuvnum={1}
\ificnum={1}
\date={January, 1996}
\titlepage
\title{New Generalized Poisson Structures}
\NAME
\FTUV
\abstract
New generalized Poisson structures are introduced by using
suitable skew-symmetric 
contravariant tensors of even order. The corresponding `Jacobi 
identities' are provided by conditions on these tensors, 
which may be understood as
cocycle conditions. 
As an example, we provide the linear generalized Poisson structures which
can be constructed on the dual spaces of simple Lie algebras.
\endpage

\chapter{Introduction}
About twenty years ago, Nambu \refmark{\Na} proposed a generalization of the 
standard classical Hamiltonian mechanics based on a three-dimensional 
`phase space' spanned by a canonical triplet of dynamical variables and on 
two `Hamiltonians'.
His approach was later discussed by Bayen and Flato \refmark{\BF} and in 
\refmark{\MS},\refmark{\Hir}.
The subject laid dormant until recently when a higher order extension of 
Nambu's approach, involving $(n-1)$ Hamiltonians, 
was proposed by Takhtajan \refmark{\Ta} 
(see \refmark{\Cha} for applications).

Another subject, closely related to Hamiltonian dynamics, is the study of 
Poisson structures (PS) (see \refmark{\Lich,\We,\BFFLS}) on a (Poisson) 
manifold $M$. A particular case of Poisson structures is that arising when they are defined 
on the duals of Lie algebras.
The class of the linear Poisson structures was considered by Lie himself 
\refmark{\Li,\Lie}, and has been further investigated recently 
\refmark{\GMP,\CIMP,\AP}.
In general,
the property which guarantees the Jacobi identity for the 
Poisson brackets (PB) of functions on a Poisson manifold may be 
expressed \refmark{\Lich,\Tu} as $[\Lambda,\Lambda]=0$ 
where $\Lambda$ is the bivector field which may be used to 
define the Poisson structure and 
$[\ ,\ ]$ is the Schouten--Nijenhuis bracket (SNB) \refmark{\Sc,\Ni}.
In the generalizations of Hamiltonian mechanics the Jacobi identity is 
replaced by a more complicated one (the `fundamental identity' in 
\refmark{\Ta}).

The aim of this paper is to introduce a new generalization of the standard 
PS. 
This will be achieved by replacing the skew-symmetric bivector $\Lambda$
defining the standard structure by appropriate even-dimensional skew-symmetric 
contravariant tensor fields $\Lambda^{(2p)}$, and by replacing the 
Jacobi identity 
by the condition which follows from $[\Lambda^{(2p)},\Lambda^{(2p)}]=0$.
In fact, the vanishing of the SNB of $\Lambda^{(2p)}$ with itself  
allows us to introduce a generalization of the Jacobi identity in a rather 
geometrical way, and provides us with a clue for the search of generalized PS. 
As a result, we differ from other approaches \refmark{\Na,\Ta}:  
all our generalized Poisson brackets (GPB) involve an 
{\it even} number of functions, whereas this number is 
arbitrary (three in \refmark{\Na}) in earlier extensions.
Since the most important question once a new Poisson structure is introduced 
is to present specific examples of it (in other words, solutions of the 
generalized Jacobi identities which must be satisfied), 
we shall exhibit, by generalizing the 
standard linear structure on the dual space $\g^*$ to a Lie algebra $\g$, 
the linear Poisson structures 
which may be defined on the duals of all simple Lie algebras.
The solution to this problem has, in fact, a cohomological component: 
the different tensors
$\Lambda^{(2p)}$ which can be introduced are related to Lie 
algebra cohomology cocycles.
We shall also discuss here the `dynamics' associated with the
GPB but shall leave a more 
detailed account of our theory and its cohomological background 
to a forthcoming publication \refmark{\APPB}.
	
\chapter{Standard Poisson structures}

Let us recall some facts concerning Poisson structures.
Let $M$ be a manifold and ${\cal F}(M)$ be the associative algebra 
of smooth functions on $M$.
\medskip
\noindent
{\bf Definition 2.1 (PB)\quad}
A {\it Poisson bracket} $\{ \cdot ,\cdot \}$ on ${\cal F}(M)$ is an
operation assigning to every pair of functions $f_1, f_2\in {\cal F}(M)$ a new
function $\{f_1, f_2\}\in {\cal F}(M)$, which is linear in $f_1$ and $f_2$ and
satisfies the following conditions:

\noindent
a) skew-symmetry 
$$\{f_1, f_2\} = - \{f_2, f_1\}\;,\eqn\pri$$
b) Leibniz rule (derivation property)
$$\{f, gh\} = g\{f, h\} + \{f, g\}h \;,\eqn\ter$$
c) Jacobi identity
$${1\over 2}{\rm Alt} \{f_1,\{f_2,f_3\}\}\equiv
\{f_1, \{f_2, f_3\}\} + \{f_2, \{f_3,f_1\}\} +
\{f_3, \{f_1,f_2\}\}= 0\;.\eqn\sec$$
The identities \pri,\sec\ are nothing but the axioms of a Lie algebra; thus
the space ${\cal F}(M)$ endowed with the PB
$\{ \cdot ,\cdot \}$
becomes an (infinite-dimensional) Lie algebra, and $M$ is a 
{\it Poisson manifold}.

Let $x^{j}$ be local coordinates on $U\subset M$ and consider PB
of the form
$$
\{f(x), g(x)\} = \omega ^{jk}(x)\partial _{j}f\partial _{k}g\quad,\quad 
\partial _{j} = {\partial \over {\partial x^{j}}}\quad,\quad
j,k=1,\ldots, {\rm dim}M 
\eqn\cua
$$
Since Leibniz's rule is automatically fulfilled, $\omega^{ij}(x)$ defines a 
PB if
$\omega ^{ij}(x) = -\omega ^{ji}(x)$
(eq. \pri) and eq. \sec\ is satisfied \ie, if
$$\omega ^{jk}\partial _{k}\omega ^{lm} + \omega ^{lk}\partial _{k}\omega 
^{mj} + \omega ^{mk}\partial _{k}\omega ^{jl} = 0\;.\eqn\sex$$

The requirements \pri\ and \ter\ imply that the PB may be given in terms of a 
skew-symmetric biderivative, \ie\ by a skew-symmetric 
bivector field (`Poisson bivector')
$\Lambda\in\wedge ^2(M)$. Locally,
$$\Lambda = {1\over 2} 
\omega ^{jk}\partial _j\wedge \partial _k\quad.\eqn\nueva
$$
Condition \sex\ may be expressed in terms of $\Lambda$ as 
$[\Lambda, \Lambda]=0$
\refmark{\Lich, \Tu}. 
A skew-symmetric tensor field $\Lambda\in\wedge^2(M)$  
such that $[\Lambda,\Lambda]=0$ defines a {\it Poisson structure} on $M$
and $M$ becomes a {\it Poisson manifold}. 
The PB is then defined by 
$$
\{f,g\}=\Lambda(df,dg)\quad,\quad f,g\in{\cal F}(M)
\quad.\eqn\oneone
$$
Two PS $\Lambda_1, \Lambda_2$ on $M$ are  
{\it compatible} if any linear combination of them is again a PS. 
In terms of the SNB this means that 
$[\Lambda_1,\Lambda_2]=0.$

Given a function $H$, the vector field $X_H=i_{dH}\Lambda$ 
(where 
$i_\alpha\Lambda(\beta):=\Lambda(\alpha,\beta)\,,\,\allowbreak\alpha,\beta$ 
one-forms), is called a
{\it Hamiltonian vector field} of $H$. 
From the 
Jacobi identity \sec\ easily follows that 
$$[X_f,X_H]=X_{\{ f,H\} }\quad.\eqn\sep$$
Thus, the Hamiltonian vector fields form a Lie subalgebra of the Lie algebra 
$\campo(M)$ of all smooth vector fields on $M.$
In local coordinates
$$X_{H} = \omega ^{jk}(x) \partial _{j}H\partial _k\quad;
\quad
X_H.f =\{H,f\}\;.\eqn\oct$$
We recall that the tensor $\omega ^{jk}(x)$ appearing in \cua, \nueva\ 
does not need to
be nondegenerate; in particular, the dimension of a Poisson manifold 
$M$ may be odd.
Only when $\Lambda$ has constant rank $2q$ (is {\it regular}) 
and the codimension
(${\rm dim}M-2q$) of the manifold is zero, $\Lambda$ defines a
{\it symplectic structure}.

\chapter{Linear Poisson structures}

A real finite--dimensional Lie algebra $\g$ with Lie bracket
$[\, .,.]$ defines in a natural way a PB
$\{\, .,.\}_\g$
on the dual space $\g^*$ of $\g$.  The natural identification $\g \cong
(\g^*)^*$, allows us to think of $\g$ as a subset of the ring of smooth
functions ${\cal F}(\g^*)$.  Choosing a linear basis
$\set{e_i}_{i=1}^r$ of $\g$, and identifying its components with linear
coordinate functions $x_i$ on the dual space $\g^*$ by means of $x_i(x) =
\langle x, e_i\rangle$ for all $x\in \g^*$, the
fundamental PB on $\g^*$ may be defined by 
$$\pois{x_i}{x_j}_\g = C_{ij}^{k} x_k \quad,\quad i,j,k=1,\ldots, r={\rm dim}
\g\quad,
\eqn\IVv$$
using that $[e_i, e_j] =  C_{ij}^{k} e_k$, where $C_{ij}^{k}$ are the
structure constants of $\g$.
Intrinsically, the PB
$\pois . . _\g$ 
on 
${\cal F}(\g^*)$ 
is defined by
$$
\pois{f}{g}_\g (x) = \langle x,[df(x),dg(x)]\rangle\quad,\quad
f,g\in {\cal F}(\g^*), x\in \g^*\quad;
\eqn\IVvi
$$
locally, 
$[df(x),dg(x)]=e_kC_{ij}^k {\partial f\over\partial x_i}
{\partial g\over \partial x_j}\,,
\, \{f,g\}_\g (x)=x_kC_{ij}^k{\partial f\over\partial x_i}
{\partial g\over \partial x_j}$.
The above PB
$\pois . . _\g$ 
is commonly called a {\it Lie--Poisson bracket}.
It is associated to the bivector field $\Lambda_\g$ on $\g^*$ 
locally written as 
$$\Lambda_\g = C_{ij}^{k} x_k \dd{x_i}\wedge \dd{x_j}
\equiv \omega_{ij}\partial^i\wedge\partial^j
\eqn\IVvii$$
(cf. \nueva), so that (cf. \oneone)
$\Lambda_\g (df\wedge dg ) = \pois{f}{g}_\g$.
It is convenient to notice here that
$[\Lambda_\g,\Lambda_\g]_S=0$ (cf. \sex) is just the Jacobi identity for $\g$,
which may be written as
$$
{1\over 2}
{\rm Alt}(C^\rho_{i_1 i_2}C^\sigma_{\rho i_3})\equiv
{1\over 2}
\epsilon^{j_1j_2j_3}_{i_1i_2i_3}C^\rho_{j_1 j_2}C^\sigma_{\rho j_3}=0\quad.
\eqn\VIa
$$

Let 
$\beta$ be a closed one form on 
$\g^*$.
The associated vector field  
$$
X_\beta= i_\beta\Lambda_\g \quad,\eqn\IVx
$$
is an infinitesimal automorphism of $\Lambda_\g$ \ie,  
$$
L_{X_\beta} \Lambda_\g = 0\quad,\eqn\IVxa
$$
and 
$[X_f, X_g] = X_{\pois f g}$ (eq. \sep);
this is proved easily using that $L_{X_f}g = \pois f g $ and
$L_{X_f}\Lambda_\g= 0$.  
It follows from \IVvii\ that
the Hamiltonian vector fields 
$X_i=i_{dx_i}\Lambda_\g$
corresponding to the linear coordinate functions $x_i$, have the
expression (cf. \oct)
$$X_i = C_{ij}^{k} x_k \dd{x_j}\quad,\quad i=1,\ldots,{\rm dim}\g\eqn\IVxii$$
so that the Poisson bivector can be written as 
$$\Lambda_\g = X_i \wedge \dd{x_i}\quad;\eqn\IVxiii$$
notice that this way of writing $\Lambda_\g$ is of course not unique.
Using the adjoint representation of 
$\g\,,\, (C_i)^k_{.j}=C^k_{ij}$
the Poisson bivector 
$\Lambda_\g$ may be rewritten as 
$$
\Lambda_\g = X_{C_i} \wedge \dd{x_i}\quad
\quad (X_{C_i}=x_k(C_i)^k_{.j}\dd{x_j})\quad;
\eqn\IVxvi
$$
the vector fields $X_{C_i}$ provide a realization of ${\rm ad}\g$
in terms of vector fields on $\g^*$.
 
\chapter{Generalized Poisson structures}
A rather stringent condition needed to define a PS on a manifold is the 
Jacobi identity \sec.
In terms of $\Lambda$, this condition is given in a convenient 
geometrical way by the vanishing of the SNB of $\Lambda\equiv\Lambda^{(2)}$
with itself, 
$[\Lambda^{(2)}, \Lambda^{(2)}]=0$.
So, it seems natural to consider generalizations of the standard PS
in terms of $2p$-ary operations determined by
skew-symmetric $2p$-vector fields $\Lambda ^{(2p)}$, the case $p=1$ 
being the standard 
one. Since the SNB of two skew-symmetric contravariant tensor fields $A, B$
of degree 
\foot{
Notice that the algebra of multivector fields is a graded superalgebra and 
that the {\it degree} of a multivector $A$ is equal to (${\it order\,A}-1$).
Thus, the standard PS defined by 
$\Lambda$ is of even order (two) but of odd degree 
(one).
} 
$a,b$ 
satisfies $[A,B]=-(-1)^{ab}[B,A]$, only $[\Lambda',\Lambda']= 0$ 
for $\Lambda'$ of odd degree will be 
meaningful since this SNB will vanish identically if $\Lambda'$ is of even
degree.

Having this in mind, let us introduce first the GPB.
\medskip
\noindent
{\bf Definition 4.1 \quad} A generalized Poisson bracket
$\{\cdot ,\cdot ,\ldots ,\cdot ,\cdot \}$ 
on $M$ is a mapping
${\cal F}(M)\times \ld^{2p}\times{\cal F}(M)\to {\cal F}(M)$
assigning a function $\{f_1, f_2,\ldots ,f_{2p}\}$ to every set 
$f_1,\ldots ,f_{2p}\in {\cal F}(M)$
which is linear in
all arguments and satisfies the following conditions:
\medskip
\noindent
a) complete skew-symmetry in $f_j$;
\medskip
\noindent
b) Leibniz rule: $\forall f_i,g,h\in {\cal F}(M),$
$$
\{f_1,f_2,\ldots ,f_{2p-1},gh\} = g\,\{f_1,f_2,\ldots ,f_{2p-1},h\}
+\{f_1,f_2,\ldots ,f_{2p-1},g\}h\quad;
\eqn\fourone
$$
\medskip
\noindent
c) generalized Jacobi identity: $\forall f_i\in {\cal F}(M),$
$$
{\rm Alt}\,\{f_1,f_2,\ldots ,f_{2p-1}\{
f_{2p},\ldots ,f_{4p-1}\}\} = 0\quad.
\eqn\fourtwo
$$

Conditions a) and b) imply that our GPB is given by a
skew-symmetric multiderivative, \ie\ by an completely skew-symmetric
$2p$-vector field $\Lambda ^{(2p)}\in \wedge^{2p}(M)$.
Condition \fourtwo\ will be called the {\it generalized
Jacobi identity}; for $p=2$ it contains $35$ terms ($C^{2p-1}_{4p-1}$ in the 
general case). It may be rewritten as $[\Lambda ^{(2p)}, \Lambda ^{(2p)}]
=0$;
$\Lambda^{(2p)}$ 
defines a GPB.
Clearly, the above relations reproduce the ordinary case \pri--\sec\ for
$p=1$. The compatibility condition of the standard case may be now extended in
the following sense:
two generalized Poisson structures 
$\Lambda ^{(2p)}$ and ${\Lambda}^{(2q)}$ on $M$
are called {\it compatible} if they `commute' \ie, $[\Lambda ^{(2p)},
{\Lambda}^{(2q)}]=0$. 
Let us emphasize that this generalized Poisson structure is different from the
Nambu structure \refmark{\Na} recently generalized in \refmark{\Ta}.
Moreover, we shall see in Sec. 5 
that our generalized linear PS are automatically
obtained from {\it constant} skew-symmetric tensors of order $2p+1$.

Let $x^j$ be local coordinates on $U\subset M$. Then the GPB
has the form
$$
\{f_1(x), f_2(x),\ldots ,f_{2p}(x)\}=\omega _{j_1j_2\ldots j_{2p}}
\partial^{j_1}f_1\,\partial^{j_2}f_2\,\ldots \,\partial^{j_{2p}}f_{2p}\;\quad.
\eqn\fourthree
$$
where 
$\omega_{j_1j_2\ldots j_{2p}}$ 
are the coordinates of a completely skew-symmetric tensor which satisfies
$$
{\rm Alt}\,(\omega _{j_1j_2\ldots j_{2p-1}k}\,\partial ^k\,
\omega_{j_{2p}\ldots j_{4p-1}}) = 0
\eqn\fourfour
$$
as a result of \fourtwo. In terms of a skew-symmetric tensor field of order 
$2p$
the generalized Poisson structure is defined by 
$$
\Lambda ^{(2p)} = {1\over {(2p)!}}\,\omega _{j_1\ldots j_{2p}}\,\partial^{j_1}
\land \ldots \land \partial ^{j_{2p}}\quad.
\eqn\fourfive
$$
Then, it is easy to check that the vanishing of the SNB
$[\Lambda ^{(2p)},\Lambda ^{(2p)}]=0$ reproduces eq. \fourfour.

Let us now define the dynamical system associated with the above generalized 
Poisson structure. 
Namely, let us fix a set of $(2p-1)$ `Hamiltonian' functions 
$H_1,H_2,\ldots , H_{2p-1}$ 
and consider the system
$$
\dot x_j=\{H_1,\ldots ,H_{2p-1}, x_j\}\quad,
$$
or, in general,
$$
\dot f=\{H_1,\ldots ,H_{2p-1}, f\}\quad.
\eqn\foursix 
$$
\medskip
\noindent
{\bf Definition 4.2\quad}
A function $f\in{\cal F}(M)$ is a constant of motion if
\foursix\ is zero.

Due to the skew-symmetry, the `Hamiltonian'
functions $H_1,\ldots ,H_{2p-1}$ are all constants of motion
but the system may have additional ones $h_{2p},\ldots ,h_k;
\,\,k\geq 2p$. 
\medskip
\noindent
{\bf Definition 4.3\quad}
A set of functions $(f_1,\ldots ,f_k)\,,k\geq 2p$ 
is in {\it involution} if the GPB vanishes for any subset of $2p$ functions.

Let us note also the following generalization of the Poisson theorem
\refmark{\Po}.
\medskip
\noindent
{\bf Theorem 4.1\quad}
Let $f_1,\ldots ,f_{q}\,,\, q\geq 2p$ be such that the
set of functions 
$(H_{1},\ldots ,\allowbreak H_{2p-1},
f_{i_{1}},\ldots ,f_{i_{2p-1}})$ is in involution 
(this implies, in particular, that the $f_i\,,\, i=1,\ldots,q$ 
are constants of motion).
Then the
quantities $\{f_{i_1},\ldots ,f_{i_{2p}}\}$ are also constants of motion.
\medskip
\noindent
{\bf Definition 4.4\quad}
A function $c(x)$ will be called a {\it Casimir function}
if $\{g_1, g_2,\ldots ,\allowbreak
g_{2p-1},\,c\}=0$ for any set of functions $(g_1, g_2,\ldots ,g_{2p-1})$.
If one of the Hamiltonians 
$(H_1,\ldots ,H_{2p-1})$ is a Casimir function, then 
the generalized dynamics defined by \foursix\ is trivial.

As an example of these generalized Poisson structures
we now show succinctly that any simple Lie algebra $\g$ of rank $l$ 
provides a family
of $l$ generalized linear Poisson structures, 
and that each of them may be characterized
by a cocycle in the Lie algebra cohomology. 

\chapter{Generalized Poisson structures on the duals of simple Lie algebras}
Let $\g$ be the Lie algebra 
of a simple compact group $G$. In this case the de Rham 
cohomology ring on the group manifold $G$ is the same as the Lie algebra 
cohomology ring $H^*_0(\g,\R)$ for the trivial action. In its
Chevalley-Eilenberg version the Lie algebra cocycles are represented by 
bi-invariant (\ie, left and right invariant and hence closed) forms on $G$
\refmark{\CE} (see also, \eg, \refmark{\AI}).
For instance, if using the Killing metric $k_{ij}$ 
we introduce the skew-symmetric order three tensor
$$
\omega(e_i,e_j,e_k)
:=k([e_i,e_j],e_k)=C_{ij}^{l}k_{lk}=C_{ijk}\,,\, 
e_i\in\g\, (i,j,k=1,\ldots,r={\rm dim}\g)
\eqn\VIi
$$
this defines by left translation a left-invariant (LI) form 
on $G$ which is also right-invariant. The bi-invariance of $\omega$ then 
reads
$$
\omega([e_l,e_i],e_j,e_k)+\omega(e_i,[e_l,e_j],e_k)+
\omega(e_i,e_j,[e_l,e_k])=0
\quad,
\eqn\VIii
$$
where the $e_i$ are now understood as LI vector fields on $G$ obtained 
by left translation from the corresponding basis of 
$\g=T_e(G)$. Eq \VIii\ (the Jacobi identity) thus implies 
a three cocycle condition on $\omega$; as a result $H^3_0(\g,\R)\neq 0$ for any 
simple Lie algebra as is well known.
In terms of the standard Poisson structure, this means that the linear
structure defined by \IVvii\ is associated with a non-trivial two-cocycle on
$\g$ and that $[\Lambda^{(2)},\Lambda^{(2)}]=0$ (eq. \VIa) is precisely
the cocycle condition. This indicates that the generalized linear Poisson 
structures on $\g^*$ may be found by looking for higher order cocycles.

The cohomology ring of any simple Lie algebra of rank $l$ is a free ring 
generated by the $l$ (primitive) forms on $G$ of odd order $(2m-1)$. 
These forms are associated 
with the $l$ primitive symmetric invariant tensors $k_{i_1\ldots i_m}$ 
of order $m$ which may be defined on 
$\g$ and of which the Killing tensor $k_{i_1i_2}$ is just the 
first example.
For the $A_l$ series $(su(l+1))$, for instance, these forms have order 
$3,5,\ldots,(2l+1)$; other orders (but always including $3$)
appear for the different simple algebras 
(see, \eg, \refmark{\AI}).
As a result, it is possible to associate a $(2m-2)$ 
skew-symmetric contravariant primitive 
tensor field linear in $x_j$ 
to each symmetric invariant polynomial $k_{i_1\ldots i_m}$ 
of order $m$.
The case $m=2$ leads to the $\Lambda^{(2)}$ of \IVvii, \IVxiii.
We shall not describe the theory in detail, 
and limit ourselves to illustrate the main theorem below with an example.
\medskip
\noindent
{\bf Theorem 5.1.\quad}
Let $\g$ be a simple compact algebra, and let 
$k_{i_1\ldots i_m}$ be
a primitive invariant symmetric polynomial of order $m$.
Then, the tensor $\omega_{\rho l_2\ldots l_{2m-2} \sigma}$
$$
\omega_{\rho l_2\ldots l_{2m-2} \sigma}:=
\epsilon^{j_2\ldots j_{2m-2}}_{l_2\ldots l_{2m-2}}
\tilde\omega_{\rho j_2\ldots j_{2m-2} \sigma}
\;,\;
\tilde\omega_{\rho j_2\ldots j_{2m-2}\sigma}:=
k_{i_1\ldots i_{m-1}\sigma}
C^{i_1}_{\rho j_2}\ldots C^{i_{m-1}}_{j_{2m-3}j_{2m-2}}
\eqn\fiveone
$$
is completely skew-symmetric, defines a Lie algebra cocycle 
\foot
{The origin of \fiveone\ is easy to understand since given a symmetric 
invariant polynomial $k_{i_1\ldots i_m}$ on $\g$, the associated 
skew-symmetric multilinear tensor $\omega_{i_1\ldots i_{2m-1}}$ 
is given by
$$
\omega(e_{i_1},\ldots, e_{i_{2m-1}})=
\sum_{s\in S_{(2m-1)}}{\pi(s)}
k([e_{s(i_1)},e_{s(i_2)}],
[e_{s(i_3)},e_{s(i_4)}],\ldots,
[e_{s(i_{2m-3})},e_{s(i_{2m-2})}],e_{s(i_{2m-1})})
$$
where $\pi(s)$ is the parity sign of the permutation $s\in S_{(2m-1)}$.}
on $\g$ and 
$$
\Lambda^{(2m-2)}={1\over (2m-2)!}{\omega_{l_1\ldots l_{2m-2}}}^\sigma
x_\sigma \partial^{l_1}\wedge\ldots\wedge \partial^{l_{2m-2}}
\eqn\fivetwo
$$
defines a generalized Poisson structure on $\g$. 

\noindent
{\it Proof:\quad} The theorem is proved using that the SNB
$[\Lambda^{(2m-2)},\Lambda^{(2m-2)}]$
is zero due to the cocycle condition satisfied by
$\omega_{\rho l_2\ldots l_{2m-2} \sigma}$. In particular 
$$
\{x_{i_1},x_{i_2},\ldots,x_{i_{2m-2}}\}=
{\omega_{i_1\ldots i_{2m-2}}}^\sigma x_\sigma
\eqn\inpart
$$
where
${\omega_{i_1\ldots i_{2m-2}}}^\sigma$ are the `structure constants' 
defining the $(2m-1)$ cocycle and hence the generalized PS.
In fact, it may be shown that different 
$\Lambda^{(2m-2)}\,,\,\Lambda^{(2m'-2)}$ tensors also commute with 
respect to the SNB and that they generate a free ring.

\noindent
{\it Note.\quad}
The requirement of compactness is introduced to have a definite 
Killing--Cartan metric 
which then may be taken as the unit matrix; this allows us to identify
upper and lower indices.
\medskip
\noindent
{\bf Example ({\it Generalized PS on} $su(3)^*$)\quad}
Let $\g=su(3)$. Besides the Killing metric
(which leads to the standard linear PS on the dual space
$su(3)^*$),
$su(3)$ admits another symmetric $ad$--invariant polynomial which may be 
expressed as ${\rm Tr}(\lambda_i\{\lambda_j,\lambda_k\})=4 d_{ijk}$
(the $d_{ijk}$ are the constants appearing in the anticommutator of the
Gell--Mann $\lambda_i$ matrices,
$\{\lambda_i,\lambda_j\}={4\over 3}\delta_{ij}1_3+2d_{ijk}\lambda_k$).
Then, the new Poisson structure is defined by
$$
\Lambda^{(4)}={1\over 4!}{\omega_{i_1i_2i_3i_4}}^\sigma x_\sigma
\dd{x_{i_1}}\wedge\ldots\wedge\dd{x_{i_4}}
\quad,\quad
\omega_{\rho i_2i_3i_4 \sigma}:={1\over 2}
\epsilon^{j_2j_3j_4}_{i_2i_3i_4}d_{k_1k_2\sigma}C^{k_1}_{\rho j_2}
C^{k_2}_{j_3j_4}\quad.\eqn\fivethree
$$
In fact, the 
$\omega_{\rho j_2j_3j_4 \sigma}$ 
in \fivethree\ is what appears
in the `four--commutators'
$$
[T_{j_1},T_{j_2},T_{j_3},T_{j_4}]= {\omega_{j_1j_2j_3j_4}}^\sigma T_\sigma
\quad (T_i={\lambda_i \over 2})
$$
which are given by the sum 
$\sum_{s\in S_4}\pi(s)T_{s(j_1)}T_{s(j_2)}T_{s(j_3)}T_{s(j_4)}$
of the $4!=24$ products of four $T$'s, each one with the sign dictated 
by the parity $\pi(s)$ of the 
permutation $s\in S_4$ and which give, as the Lie algebra commutator does, an 
element of $\g$ in the right-hand side.
It is not difficult now to check, using the symmetry of the $d$'s and the 
properties of the structure constants (including the Jacobi identity) that 
$[\Lambda^{(4)},\Lambda^{(4)}]=0$. 
Thus, all properties of Def. 4.1 are fulfilled and 
$\Lambda^{(4)}$ 
defines a GPB.
We refer 
to \refmark{\APPB} for further details
concerning the mathematical structure of the GPB and the contents
of the associated generalized dynamics and its quantization. We shall 
conclude here by saying that this analysis could be extended to
Lie superalgebras and super-Poisson structures.

\bigskip
\noindent
{\bf Acknowledgements.\quad}
This research has been partially supported by the CICYT (Spain) under grant 
AEN93--187. 
A. P. wishes to thank the Vicerectorate of research of Valencia University for 
making his stay in Valencia possible;
J.C.P.B. wishes to acknowledge a FPI grant from the Spanish Ministry of 
Education and Science and the CSIC.


\refout

\end